\documentclass[12pt,reqno]{amsart}

\usepackage[cp1251]{inputenc}
\usepackage[T2A]{fontenc}
\usepackage[english]{babel}

\usepackage{amsmath,amsxtra,amssymb,eufrak}
\usepackage[all]{xy}
\usepackage{mathrsfs}
\usepackage{paralist}
\usepackage[active]{srcltx} 

\theoremstyle{plain}
\newtheorem{theorem}{Theorem}
\newtheorem{lemma}[theorem]{Lemma}
\newtheorem{proposition}[theorem]{Proposition}

\theoremstyle{definition}

\newtheorem{remark}[theorem]{Remark}

\newtheorem{question}[theorem]{Question}

\begin{document}

\author{O.\,Yu.~Aristov}
\email{aristovoyu@inbox.ru}

\title{The relation ``commutator equals function'' in Banach algebras}

\markboth{O.\,Yu.~Aristov}{The relation ``commutator equals function''}

 \maketitle

\begin{abstract}
The relation $xy-yx=h(y)$, where $h$ is a holomorphic function, occurs naturally in
the definitions of some quantum groups. To attach a rigorous meaning to the right-hand side of
this equality, we assume that $x$ and $y$ are elements of a Banach algebra (or of an Arens--Michael
algebra). We prove that the universal algebra generated by a commutation relation of this kind can
be represented explicitly as an analytic Ore extension. An analysis of the structure of the algebra
shows that the set of holomorphic functions of $y$ degenerates, but at each zero of $h$, some local
algebra of power series remains. Moreover, this local algebra depends only on the order of the zero.
As an application, we prove a result about closed subalgebras of holomorphically finitely generated
algebras.
\end{abstract}

In quantum algebra, there arise commutation relations involving not only polynomials but also
more general holomorphic functions. Similar relationships can be obtained by deforming the universal
enveloping algebra of a semisimple Lie algebra (following Drinfeld and Jimbo; see the monograph~\cite[Definition~17.2.3]{Ka91})  or of the Lie algebra of the group of affine transformations of the line; see the paper~\cite{AiSa} by Aizawa and Sato. In the latter case, the condition has the simplest form:
\begin{equation}\label{relat}
[x,y]=h(y) ,
 \end{equation}
where $h$ is the hyperbolic sine \cite[formula~(3.1)]{AiSa}. (Here $[x,y]\!:=xy-yx$.)
We consider a general relation of this form, assuming that $h$ is a function holomorphic in some domain.

Certainly, such equalities make no sense for arbitrary algebras\footnote{Here and in what follows, by ``algebra'' we mean ``associative algebra with unit over the field~$\mathbb{C}$ of complex numbers.''}. To overcome this difficulty, quantum
algebras over the ring of formal series in the quantization parameter are traditionally introduced. The
alternative analytical approach of considering these relations in algebras for which the holomorphic
functional calculus theorem ensures the existence of $h(y)$ (in particular, in Banach algebras) seems to
be more natural. This point of view supposes a specialization of the quantization parameter to a complex
number. Further, one must consider the universal algebra generated by those elements $x$ and $y$ for which
the required relation holds, assuming additionally that these elements are contained in a Banach algebra
with unit and $h$ is holomorphic in some neighborhood of the spectrum of $y$.
Before investigating the properties of the resulting ``analytical form'' of the quantum group, it is necessary to study the question of whether this form is nontrivial and how rich is its structure.

The main objective of this paper is to answer this question by providing a complete description of
the universal algebra generated by elements $x$ and $y$ satisfying~\eqref{relat} for an arbitrary function $h$. Such
problems do not always have a solution in the class of Banach algebras without additional conditions on
the norms or the spectrum of elements; however, a solution can be found among algebras approximated
by Banach algebras, namely, among the Arens--Michael algebras.  (Recall that a complete locally convex
algebra is called an \emph{Arens--Michael algebra}  if its topology is generated by a family of submultiplicative
seminorms, i.e., of seminorms $\|\cdot\|$ satisfying the condition
for any elements $a$ and $b$.)

Certainly, particular cases of~\eqref{relat} have been studied from different points of view, at least for
polynomials.  For example, it has been established that if $h$ is a nonconstant polynomial, then the problem
of classification (up to similarity) of pairs of finite-dimensional operators satisfying this condition is very
difficult (it is wild according to Ostrovskii and Samoilenko \cite[Theorem~1]{OS13}).
Some spectral properties of a pair of elements of a Banach algebra satisfying~\eqref{relat} were studied by Turovskii \cite{Tu85} (cited in \cite{Shu94})); see also Shulman's survey \cite[Sec.~3]{Shu94}. However, as far as the author knows, the universal Arens--Michael algebra
generated by this relation has previously been identified only for constant or linear functions $h$. Recall
that a classical result states that if $h$ is a nonzero constant, then a nontrivial realization of relation~\eqref{relat} in a Banach algebra is impossible. The case of $h(y)=y$ was treated by Pirkovski \cite{Pir_qfree}.

Although \eqref{relat} does not imply upper bounds for the norms of powers of $x$, it rather rigidly determines
the asymptotic behavior of  $\|(y-\lambda)^n\|$ as $n\to \infty$, where $\lambda$ is some zero of the function $h$.
For example, it follows from  $[x,y]=y$ that
$$
 n\|y^n\|\le 2\|x\|\,\|y^n\|\qquad\text{for all
$n\in\mathbb{N}$}
$$
 and thus $\|y^n\|=0$ beginning with some $n$; cf.~\cite[Example 5.1]{Pir_qfree}.  This observation was used in
\cite[Proposition~5.2]{Pir_qfree} in describing the corresponding universal algebra.

The relation $[x,y]=y^2$ also implies a restriction on growth: arguing by induction, we can show that
$$
\|y^n\|= O\biggl(\frac{\|x\|^n}{n!}\biggr) \qquad (n\to \infty).
$$
A classical example is a pair of operators on $L^2[0,1]$, the operator $T$ of multiplication by the independent
variable and the indefinite integration operator
$$
 Vf(x)=\int_0^x  f(t)\, dt
$$
(which is a special case of the Volterra operator). It can readily be seen that $[T, V] = V^2$. There is a vast
literature, devoted to calculation and estimation of the norms of powers of $V$; of special mention are the
papers
\cite{LW97}--\cite{LR98}. In particular, it has been proved that
$\lim_{n\to\infty} n!\,\|V^n\| = 1/ 2$; see  \cite[Theorem~5.4]{Ke99} and
\cite[Remark~3]{Th98}.

A description of a universal algebra in the general case is given in Theorem~\ref{comfunenv} (see below).  As expected
from the above examples, it depends only on the zeros of $h$ and their orders. Moreover, to every zero, there
corresponds a subalgebra that consists of (not necessarily convergent) power series and is local.

Note also that~\eqref{relat}  admits the following natural generalization.
\begin{question}
Let $U$  be a domain in $\mathbb{C}$, let $h$  be a nonzero function holomorphic on $U$, and let $\alpha$  be a
continuous endomorphism of the algebra of functions holomorphic on~$U$. What is the Arens--Michael
universal algebra for the relation
$$xy-\alpha(y)x=h(y)\,?$$
\end{question}

This question remains open in the general case. The special case of the quantum Weyl algebra with $h=1$ was described by Pirkovskii~\cite[Corollary~5.19]{Pir_qfree}.

Theorem~\ref{comfunenv}  has an application to the theory of holomorphically finitely generated (or, briefly, HFG)
algebras.  This class of algebras was first considered by Pirkovsky (see \cite{Pi14} and \cite{Pi3}) and is of interest
from the point of view of noncommutative geometry, since the commutative HFG algebras are Stein
algebras, i.e., algebras of holomorphic functions, or, more precisely, of global sections of the structure
sheaf on some Stein space (in the case of finite dimension of the embedding).

At present, the study of noncommutative HFG algebras remains in an embryonic state. In particular,
it is not clear how wide is the class of their closed subalgebras. However, at least, it is known that a
closed subalgebra of a Stein algebra need not be a Stein algebra, which means that the class of HFG
algebras is not stable with respect to the passage to closed subalgebras either.

When studying our main problem, a family of local algebras ${\mathscr A}_s$, $s\in[0,\infty]$, of power series arises.
We show below that, for rational values of  $s$,  such an algebra can be embedded as a closed subalgebra
in some HFG algebra (Theorem~\ref{Asclsa}), although it is not HFG for $s\ne0$ (Proposition~\ref{nonHFG}).
The following
questions remain open.
\begin{question}\label{questionsHFG}
Is ${\mathscr A}_s$ a closed subalgebra of an HFG algebra for an arbitrary irrational positive~$s$?
\end{question}

\begin{question}
\label{questionsSTsub}
Is ${\mathscr A}_s$ a closed subalgebra of a Stein algebra for
$s\in(0,+\infty]$?
\end{question}


\subsection*{The universal algebra and formulation of the main theorem}

Let us begin with the definition of a family of local algebras that plays an important role in our
reasoning. For every $s\ge0$, consider the following completion of the algebra of polynomials in a formal
variable $y$:
\begin{equation}\label{Asdefin}
{\mathscr A}_s\!:=\Bigl\{a=\sum_{n=0}^\infty  a_n y^n\! :
\|a\|_{r,s}\!:=\sum_{n=0}^\infty |a_n|\frac{r^n}{n!^s}<\infty
\;\forall r>0\Bigr\}.
\end{equation}
It can readily be seen that the restrictions of the seminorms $\|\cdot\|_{r,s}$ to the algebra of polynomials in $y$ are
submultiplicative (cf. the proof of part~(A)
of Proposition~\ref{cOUdernA} below), and thus ${\mathscr A}_s$ is an Arens--Michael
algebra. Let us denote the algebra  $\mathbb{C}[[y]]$  of all formal power series in $y$ by ${\mathscr A}_\infty$. This is also an
Arens--Michael algebra with respect to the topology generated by the system of submultiplicative
seminorms
$$
\|a\|_{m,\infty}\!:=\sum_{n=0}^m |a_n|, \qquad (m\in\mathbb{Z}_+)\,.
$$

For a given domain $U$ in $\mathbb{C}$, we denote the algebra of all holomorphic functions on~$U$ by $\mathcal{O}(U)$ and choose a nonzero
$h\in\mathcal{O}(U)$. Let $\{\lambda_j:\,j\in J\}$  be the set of all zeros of the function $h$ (without
repetitions), and let  $s_j\!:=1/(k_j-1)$, where $k_j$ stands for the order of~$\lambda_j$.   For every
$j$, we take the
algebra ${\mathscr A}_{s_j}$ and set
\begin{equation}\label{Aaspr}
{\mathscr A}\!:=\prod_{j\in J} {\mathscr A}_{s_j}.
\end{equation}

Below we show that the desired universal algebra is an analytic Ore extension of the algebra~${\mathscr A}$. This analytic version of a classical notion was suggested by Pirkovskii in~\cite[Sec.~4.1]{Pir_qfree}. Let us recall the
necessary definitions and facts. Let $d$ be a derivation of some algebra $R$.  A seminorm
$\|\cdot\|$ on~$R$  is said to
be  $d$-\emph{stable} if there is a $C>0$  such that $\|d(r)\|\le C\,\|r\|$ for all $r\in
R$ \cite[Definition~4.1]{Pir_qfree}. A derivation $d$
of an Arens--Michael algebra $R$ is said to be $m$-\emph{localizable} if the topology on $R$ is generated by a family
of $d$-stable submultiplicative seminorms \cite[Definition~4.4]{Pir_qfree}.

\begin{proposition}\label{uniprOre}{\rm \cite[Propositions~4.4 and~4.6 and Remark~4.6]{Pir_qfree}}
 Let $R$ be an Arens--Michael algebra,
and let $d$ be its $m$-localizable derivation. Then there exist an Arens--Michael algebra $E$, an
$x\in E$,  and a continuous homomorphism $\eta\!:R\to E$ such that
$[x,\eta(r)]=\eta d(r)$  for all $r\in R$  which
has the following universal property. If $B$ is an Arens--Michael algebra,  $\breve x\in B$, and $\nu\!:R\to B$
is a continuous homomorphism such that $[\breve x,\nu(r)]=\nu d(r)$ for all $r\in R$, then there is a unique
continuous homomorphism $\tau\!:E\to B$ such that
$\nu=\tau \eta$ and $\tau(x)=\breve x$.
\end{proposition}

The universal algebra $E$ is denoted by $\mathcal{O}(\mathbb{C},R;d)$ and called an
\emph{analytic Ore extension}  of the
algebra~$R$
\cite[Definition~4.3]{Pir_qfree}. (Note that this notion still makes sense in the more general case in
which $d$ is an  $\alpha$-derivation for some endomorphism~$\alpha$  of the algebra~$R$.)  The underlying locally convex
space of the algebra $\mathcal{O}(\mathbb{C},R;d)$   is the complete projective tensor product $R{\mathbin{\widehat{\otimes}}}\mathcal{O}(\mathbb{C})$,
the homomorphism $\eta\!:R\to\mathcal{O}(\mathbb{C},R;d)$   is defined by the condition
$r\mapsto r\otimes 1$, and $x=1\otimes z$, where $z$ is the identity function
on $\mathbb{C}$ \cite[Proposition~4.3]{Pir_qfree}.

We denote by~$y_j$ the corresponding formal variable in~${\mathscr A}_{s_j}$  and by $y$ the sequence $(y_j+\lambda_j:\!j\in J)$ in~${\mathscr A}$.
 It can readily be proved that  $y_j\in{\mathscr A}_{s_j}$ is quasinilpotent (since  $s_j>0$).
Thus, each of the algebras ${\mathscr A}_{s_j}$  is local, and hence the spectrum of $y$ coincides with $\{\lambda_j:\!j\in J\}$.  Since the spectrum is
contained in $U$ and ${\mathscr A}$  is an Arens--Michael algebra, it follows that the holomorphic functional calculus
for $y$ is well defined  \cite[Chapter~VI, Theorem~3.2]{Mal};  in particular, there is an  $h(y)\in {\mathscr A}$.

It can readily be seen that
\begin{equation}\label{dejdef}
\delta_j\!:=h(y_j+\lambda_j)\frac{d}{dy_j}
\end{equation}
is a derivation of~${\mathscr A}_{s_j}$, i.e.,
 $$
 \delta_j(ab)=a \delta_j(b)+\delta_j(a)b\qquad\text{for any $a,b\in {\mathscr A}_{s_j}$}.
$$
Moreover, the map
$\delta\!:{\mathscr A}\to {\mathscr A}$, being the product of all~$\delta_j$, is also a derivation.

We are now ready to state the main result of the paper, which is the following theorem on the universal
algebra.

 \begin{theorem}\label{comfunenv}
Suppose that $U$ is a domain in $\mathbb{C}$ and $h$ is a nonzero holomorphic function on U. Let
$\{\lambda_j:\,j\in J\}$
be the set of all distinct zeros of $h$, and let $s_j\!:=1/(k_j-1)$, where $k_j$ is the order of
$\lambda_j$.

\emph{(A)}~Let the Arens--Michael algebra  ${\mathscr A}$, its derivation $\delta$ and $y\in{\mathscr A}$ be defined as above. Then $\delta$ is
$m$-localizable, and hence $\mathcal{O}(\mathbb{C},{\mathscr A};\delta)$  is well defined and is an Arens--Michael algebra; moreover,
the spectrum of $y$ is contained in $U$ and the relation $[x,y]=h(y)$ holds.

\emph{(B)}~If an Arens--Michael algebra (in particular, a Banach algebra)~$B$
contains elements~$\breve x$ and~$\breve y$ such that the spectrum of~$\breve y$ is contained in~$U$ and
$[\breve  x,\breve y]=h(\breve y)$, then there is a unique continuous
homomorphism
$$
\tau\!:\mathcal{O}(\mathbb{C},{\mathscr A};\delta)\to B\,,
$$
for which $\tau(x)=\breve x$ and $\tau(y)=\breve y$.
\end{theorem}

In particular, if $h$ has no zeros, then the universal algebra is isomorphic to $\{0\}$. Thus, the theorem
includes the classical result on the relation $[x,y]=1$ as a special case. On the other hand, deforming the
function  $h(y)=y$  into $h(y)={\sinh \hbar y}/{\sinh
\hbar}$  as in~\cite{AiSa}  (here the quantization parameter is
$\hbar\in\mathbb{C}\setminus\{0\}$), we see that the universal algebra contains infinitely many copies of $\mathbb{C}[[y]]$. This follows from the
theorem, because the set of zeros of the hyperbolic sine is infinite and all zeros are of order~$1$. Since
the universal algebra in the classical case contains only one such copy, this visually demonstrates the
effect of quantization.

\begin{remark}
 If $h$ is identically equal to zero in a domain $U$, then, obviously, the universal algebra also
exists and is topologically isomorphic to
$\mathcal{O}(U){\mathbin{\widehat{\otimes}}}\mathcal{O}(\mathbb{C})$.   For the case in which $U=\mathbb{C}$, this well agrees
with our notation, since $\mathcal{O}(\mathbb{C})\cong{\mathscr A}_0$ (if we assume that we have here a ``zero of infinite order'').
\end{remark}

We will first prove part (A) of Theorem~\ref{comfunenv} and then part (B). The idea of the proof of part (A) is to
construct a family of $\delta_0$-stable submultiplicative seminorms on
$\mathcal{O}(U)$,  where the derivation $\delta_0$ is given
by the formula
$\delta_0(f)=hf'$. To establish the validity of part~(B), we will prove that this family is sufficient
for describing the topology of the universal algebra.

\subsection*{Proof of part~(A) of Theorem~\ref{comfunenv}}
On $\mathcal{O}(U)$,  we consider the seminorms (cf. the definition of the algebras ${\mathscr A}_s$ in~\eqref{Asdefin}):
\begin{equation}\label{norrs}
\|f\|_{\lambda,r,s}\!:=\sum_{n=0}^\infty
|f^{(n)}(\lambda)|\,\frac{r^n}{n!^{s+1}}\,,\qquad
 \|f\|_{\lambda,m,\infty}\!:=
\sum_{n=0}^m\frac{|f^{(n)}(\lambda)|}{n!}
\end{equation}
(here $\lambda\in U$, $r>0$, $s> 0$ and
$m\in\mathbb{Z}_+$) and the standard seminorms
$$
|f|_K\!:=\sup\{|f(z)|:\, z\in K\}
$$
(here $K$ stands for a compact subset of $U$).

The proofs of both parts of Theorem~\ref{comfunenv} use the following two assertions.

\begin{proposition}\label{cOUdernA}
Consider the derivation of the algebra  $\mathcal{O}(U)$ given by the formula $\delta_0(f)=hf'$.

\emph{(A)}~Each of the seminorms $\|\cdot\|_{\lambda,r,s}$ and
$\|\cdot\|_{\lambda,m,\infty}$ ($\lambda\in U$, $r>0$, $s>0$, and $m\in\mathbb{Z}_+$)
is submultiplicative and, up to a constant, is dominated by the seminorm $|\cdot|_D$  for some closed disk $D\subset
U$ of sufficiently small radius centered at~$\lambda$; hence these seminorms are continuous.

\emph{(B)}~If $\lambda$ is a zero of the function $h$, then
$\|\cdot\|_{\lambda,m,\infty}$ is $\delta_0$-stable for all
$m\in\mathbb{Z}_+$.

\emph{(C)}~If $\lambda$ is a zero of order $k>1$  of the function $h$, then $\|\cdot\|_{\lambda,r,s}$ is $\delta_0$-stable for all
$r>0$ and $s\ge1/(k-1)$.
\end{proposition}
\begin{proof}
 Let us prove (A). If s is finite, then the submultiplicativity of every seminorm of the form $\|\cdot\|_{\lambda,r,s}$
follows from the Leibniz formula and the inequality
$$
\frac{r^{l+n}}{(l+n)!^{s}}\le \frac{r^l}{l!^s}\,\frac{r^n}{n!^s}\qquad (l,n\in\mathbb{Z}_+).
$$
In the case of $\|\cdot\|_{\lambda,m,\infty}$,  the proof of submultiplicativity is straightforward.

The second part of the assertion follows from the Cauchy inequalities for the coefficients of the Taylor
series and from the fact that the topology on $\mathcal{O}(U)$  is generated by the family \{$|\cdot|_K\}$
of seminorms, where
$K$ ranges over all compact subsets of $U$.

Let us prove (B). Suppose that $h(\lambda)=0$ and choose an $m\in\mathbb{Z}_+$. Since
\begin{equation}\label{hfpr}
\delta_0(f)^{(m)}=(hf')^{(m)}=\sum_{p=0}^m {m\choose p}
h^{(p)}f^{(m-p+1)}\qquad(f\in \mathcal{O}(U)),
\end{equation}
it follows that the number $\delta_0(f)^{(m)}(\lambda)$ is a linear function in
$f(\lambda),f'(\lambda),\ldots, f^{(m)}(\lambda)$ with coefficients
independent of $f$. This readily implies that $f$
$\|\cdot\|_{\lambda,m,\infty}$ is $\delta_0$-stable for every 
$m\in\mathbb{Z}_+$.

Let us prove (C). Suppose that $\lambda$ is a zero of order $k>1$ of the function $h$. Choose an $r>0$ and $s\ge1/(k-1)$.
By the Cauchy inequalities, there are positive numbers $M$ and $R$ such that
 $|h^{(p)}(\lambda)/p!|\le
M/R^p$ for all $p\in\mathbb{Z}_+$.  Substituting \eqref{hfpr} into
\eqref{norrs}  and taking into account the fact that
$$
h(\lambda)=h'(\lambda)=\cdots=h^{(k-1)}(\lambda)=0,
$$
   we obtain
$$
\|\delta_0(f)\|_{\lambda,r,s}=\sum_{n=k}^\infty \left|\sum_{p=k}^n
{n\choose p}
h^{(p)}(\lambda)f^{(n-p+1)}(\lambda)\right|\frac{r^n}{n!^{s+1}}\le
\sum_{n=k}^\infty \sum_{p=k}^n\frac{M}{R^p}\,
\frac{|f^{(n-p+1)}(\lambda)|}{(n-p)!}\, \frac{r^n}{n!^{s}}
$$
for every $f\in \mathcal{O}(U)$.
Making the changes $p=m+k$ and $n=q+m+k-1$, we can write
\begin{multline*}
\sum_{q=1}^\infty \sum_{m=0}^\infty \frac{M}{R^{m+k}}\,
\frac{r^{q+m+k-1}}{(q+m+k-1)!^{s}}\,
\frac{|f^{(q)}(\lambda)|}{(q-1)!}=\\
=
\sum_{q=1}^\infty\frac{q\,r^{q+k-1}}{(q+k-1)!^{s}}\,\frac{|f^{(q)}(\lambda)|}{q!}\sum_{m=0}^\infty
\frac{M}{R^{m+k}} \,\frac{r^{m}(q+k-1)!^{s}}{(q+m+k-1)!^{s}}\,.
\end{multline*}
Applying the inequality
$$
\frac{(q+k-1)!^{s}}{{R^{m+k}} (q+m+k-1)!^{s}}\le \frac
1{R^{k}m!^s}
$$
to the terms of the sum over $m$, we arrive at
$$
\|\delta_0(f)\|_{\lambda,r,s}\le\frac{C_r M}{R^{k}}
\sum_{q=1}^\infty\frac{q\,
r^{q+k-1}}{(q+k-1)!^{s}}\,\frac{|f^{(q)}(\lambda)|}{q!}\,,
$$
where $C_r\!:=\sum_{m=0}^\infty r^{m}/m!^{s}$ (note that
$C_r<\infty$, because $s>0$).  Further, the condition $s\ge1/(k-1)$ implies the inequality
$$
q!\,q^{1/s}\le(q+k-1)!\,.
$$
Thus, we have obtained the final bound
$$
\|\delta_0(f)\|_{\lambda,r,s}\le\frac{C_r M}{R^{k}}\, \sum_{q=1}^\infty
|f^{(q)}(\lambda)|\, \frac{r^{q+k-1}}{q!^{s+1}}=\frac{C_r
Mr^{k-1}}{R^{k}}\,\|f\|_{\lambda,r,s}.
$$
This completes the proof of the proposition.
\end{proof}

Further, we consider the following family of seminorms on  $\mathcal{O}(U)$:
\begin{equation}\label{prsys}
\{\|\cdot\|_{\lambda_j,r_j,s_j}\!: j\in I\},
\end{equation}
where, as above, $\{\lambda_j:\,j\in J\}$  is the set of all zeros of the function $h$, $k_j$ is the order of~$\lambda_j$,
$s_j\!:=1/(k_j-1)$, and the domain of variation of the parameters is given by the following rule:  if $\lambda_j>1$, then $r_j\in\mathbb{R}_+$, and if
$\lambda_j=1$, then $s_j=\infty$ and $r_j\in\mathbb{N}$.

\begin{proposition}\label{compllem}
The homomorphism
$$
\mu\!:\mathcal{O}(U)\to {\mathscr A}\!: f\mapsto \left(\sum_{n=0}^\infty \frac{f^{(n)}(\lambda_j)}{n!}\,y_j^n\right)_j
$$
is well defined and continuous, and this is the completion homomorphism with respect to~\eqref{prsys}.
Moreover, $\mu\delta_0=\delta\mu$.
\end{proposition}
\begin{proof}
The definition of ${\mathscr A}$ (see \eqref{Asdefin} and \eqref{Aaspr})  implies that~$\mu$ is well defined and continuous and that the
topologies coincide. It remains to prove that every $(\alpha_j:\,j\in J)\in {\mathscr A}$ can be approximated by elements of
the form $\mu(f)$, where $f\in\mathcal{O}(U)$. Moreover, we may assume that  $(\alpha_j)$  is a finite sequence of polynomials.

Let us choose an $r_j$ for every $j$ such that $\alpha_j\ne0$. It follows from the finiteness of the set of these~$j$  and from part (A) of Proposition~\ref{cOUdernA} that there are a $C>0$ and a compact subset $K$ of $U$ which is a
union of closed disks
$D_j$ centered at~$\lambda_j$  such that $\|f\|_{\lambda_j,r_j,s_j} \le
C\,\|f\|_K$ for all $f\in\mathcal{O}(U)$.  Reducing the
radii if necessary, we may assume that the disks~$D_j$  are pairwise disjoint. We define a function $g$ on $K$
as follows: for every $z\in D_j$, $g$ coincides with the polynomial~$\alpha_j$ in which the substitution $y_j=z-\lambda_j$ is made. Since $K$ is compact, $\mathbb{C}\setminus K$ is connected, and $g$ is continuous on $K$ and holomorphic on its
interior, it follows that we can apply Mergelyan’s theorem, which claims that $g$ can be approximated by
polynomials uniformly on $K$ (see, e.g.,
\cite[Theorem~20.5]{Ru87}). Hence $g$ is approximated by polynomials
with respect to the topology given by~\eqref{prsys}.

The equality $\mu\delta_0=\delta\mu$ is obtained by a straightforward calculation using~\eqref{hfpr}.
\end{proof}

We can now complete the proof of the first part of the main theorem.

\begin{proof}[Proof of Theorem~\ref{comfunenv}.  Part (A)]
 Recall that $\delta$ is the product of all derivations~$\delta_j$ defined in~\eqref{dejdef}. Therefore, to prove that~$\delta$  is $m$-localizable, it suffices to prove this for every~$\delta_j$. By Proposition~\ref{compllem}  and parts (B)
and (C) of Proposition~\ref{cOUdernA}, every seminorm in the defining system for~${\mathscr A}_{s_j}$  is $\delta_j$-stable, as desired.

The fact that the spectrum of $y$ consists of those zeros of the function $h$ that belong to $U$, and hence
 $h(y)$ is well defined, was already noted above. The relation $[x,y]=h(y)$ follows from the construction of
the algebra $\mathcal{O}(\mathbb{C},{\mathscr A};\delta)$  (see Proposition~\ref{uniprOre}).
\end{proof}

\subsection*{Proof of Part (B) of Theorem~\ref{comfunenv}}
We set $\delta_0(f)=hf'$, as above. We need the following assertion.

\begin{proposition}\label{cOUdernB}
Every continuous submultiplicative $\delta_0$-stable seminorm on $\mathcal{O}(U)$  is dominated by
the family~\eqref{prsys}  of seminorms.
\end{proposition}

To prove this proposition, we use two lemmas. Below we denote by ${\mathop{\mathrm{Sp}}\nolimits}_B b$  the spectrum of the element $b$ of the algebra $B$.
\begin{lemma}\label{lem1}
 Let $d$ be a continuous derivation of a commutative Banach algebra $B$, let $b\in B$,  and $h\in\mathcal{O}(V)$, where $V$ is a domain in $\mathbb{C}$
containing ${\mathop{\mathrm{Sp}}\nolimits}_B b$. If
$d(b)=h(b)$, then ${\mathop{\mathrm{Sp}}\nolimits}_B b$  consists of
zeros of $h$ and is finite.
\end{lemma}
\begin{proof}
The Singer–Wermer theorem \cite[Theorem~7.2.10]{X2}  claims that the range of a continuous
derivation of a commutative Banach algebra is contained in the Jacobson radical. In particular, it follows
from $h(b)=d(b)$ that $h(b)$  belongs to the radical. Since $B$ is commutative, every element of the radical
is topologically nilpotent \cite[Theorem~2.1.34]{X2};
hence the spectrum ${\mathop{\mathrm{Sp}}\nolimits}_B h(b)$ is~$\{0\}$.
Since ${\mathop{\mathrm{Sp}}\nolimits}_B b\subset V$, it follows from the spectral mapping theorem
\cite[Theorem~2.2.23]{X2}  that ${\mathop{\mathrm{Sp}}\nolimits}_B b$
is contained in the set of
zeros of $h$. Since this set does not have limit points in $V$ and ${\mathop{\mathrm{Sp}}\nolimits}_B b$  is compact, it is finite.
\end{proof}

Let $\|\cdot\|$
be a continuous submultiplicative $\delta_0$-stable seminorm on
$\mathcal{O}(U)$. We denote the completion
of   $\mathcal{O}(U)$ with respect to $\|\cdot\|$
by $B$ and the image of the identity function under the completion
homomorphism $\mathcal{O}(U)\to B$
by~$\breve y$. Obviously, $\delta_0$ extends to a continuous derivation of $B$. We denote
the norm and the derivation extended to $B$ by the same symbols.

Choose a $\lambda\in {\mathop{\mathrm{Sp}}\nolimits}_B \breve y$.  By Lemma~\ref{lem1},
the number $\lambda$  is a zero of~$h$. Let $k$ denote the order of this zero.
Then there is a $g\in\mathcal{O}(U)$ such that
 $$g(\lambda)\ne
0\quad\text{and}\quad h(z)=(z-\lambda)^{k}g(z)$$ for all $z\in U$.
Let $V$ be an open neighborhood of $\lambda$ contained in $U$ and such that $g(z)\ne 0$ for $z\in V$. Then
$V$ contains no other points of  ${\mathop{\mathrm{Sp}}\nolimits}_B \breve y$. Since ${\mathop{\mathrm{Sp}}\nolimits}_B \breve y$  is finite, it follows from the holomorphic functional
calculus theorem that the characteristic function $\chi$ of the set $V$ can be applied to $\breve y$. For brevity, we will
write $g$ and $\chi$
instead of $g(\breve y)$ and $\chi(\breve y)$, respectively.

\begin{lemma}\label{lem2}

\emph{(A)}~If $k=1$, then  $\|(y-\lambda)^n \chi\|=0$ for all sufficiently large~$n$.

\emph{(B)}~If $k\ge 2$, then there are $K>0$ and $r>0$ such that
\begin{equation}\label{inkg2}
\|(y-\lambda)^n \chi\|\le K \,\frac{r^n}{\sqrt[k-1]{(n+k-1)!}}\qquad
(n\ge k-1)\,.
\end{equation}
\end{lemma}
\begin{proof}
 Without loss of generality, we can assume that $\lambda=0$.
Let us first estimate the norm of  $\breve y^n g\chi$. Let $C>0$ be such that
$\|\delta_0(b)\|\le C\|b\|$  for all $b\in B$. Further, we note that, first,  $\chi$
is an idempotent and,
therefore, it follows from $\delta_0(\chi)=2\delta_0(\chi)\chi$ and
$(1-2\chi)^2=1$ that $\delta_0(\chi)=0$. Thus,
$\delta_0(b\chi)=\delta_0(b)\chi$ for
every $b\in B$. Second, since the function $g$ has no zeros in $V$, there is a $w\in B$ such that $g\chi w=\chi$. Therefore, for every  $n\in\mathbb{N}$,  we obtain the following equality from $h(z)=z^{k}g(z)$:
$$
\delta_0(\breve y^n \chi)=\delta_0(\breve y^n) \chi=n \breve y^{n+k-1} g\chi.
$$
Hence
\begin{equation}\label{inq1}
\|\breve y^{n+k-1} g\chi\|\le n^{-1}C\,\|\breve y^n  \chi\|\le n^{-1}C\,\|\breve y^n
g\chi  \| \, \|w\|\,.
\end{equation}

Let us prove (A). Suppose that $k=1$. Then $n+k-1=n$  and, for sufficiently large $n$, we obtain
 $\|y^n g\chi\|=0$, whence $\|y^n
\chi\|\le\|y^n g\chi \| \, \|w\| =0$.

Let us prove (B). Suppose that $k\ge 2$. Let $m\in\{0,\ldots,k-1\}$, and let
$j\in\mathbb{N}$. Applying inequality \eqref{inq1} $j$ times, we obtain
$$
\|\breve y^{m+j(k-1)}\chi\|\le \|\breve y^{m+j(k-1)}g\chi\|\,\|w\| \le
\frac{C^j\,\|w\|^{j+1}\,\|\breve y^m
g\chi\|}{m(m+k-1)\cdots(m+(j-1)(k-1))}\,.
$$
Every positive integer $n$ which is not less than $k-1$ can be written in the form $m+j(k-1)$ with the
above conditions on $m$ and $j$. Since
$$m^{k-1}(m+k-1)^{k-1}\cdots(m+(j-1)(k-1))^{k-1}\ge
(m+(j-1)(k-1))!\,,$$
there are $K>0$ and  $r>0$ for which~\eqref{inkg2} holds.
\end{proof}

\begin{proof}[Proof of Proposition~\ref{cOUdernB}]
We use the notation introduced before Lemma~\ref{lem2}.
 It follows from Lemma~\ref{lem1} that ${\mathop{\mathrm{Sp}}\nolimits}_B \breve y=\{\lambda_1,\ldots,\lambda_l\}$,
where $\lambda_1,\ldots,\lambda_l$ are pairwise distinct zeros of $h$. For every  $j\in\{1,\ldots,l\}$ there is a function $g_j$ in $\mathcal{O}(U)$
 such that
 $$g_j(\lambda_j)\ne 0\quad\text{and}\quad h(z)=(z-\lambda_j)^{k_j}g_j(z)\,,$$ where
$k_j$
 is the order of $\lambda_j$  as a zero of $h$. Let us choose, for every $j$, a neighborhood $V_j$ of $\lambda_j$ such that $g_j(z)\ne 0$ in~$V_j$ and
$V_j\subset U$. It can be assumed that the neighborhoods $V_1,\ldots,V_l$ are pairwise
disjoint. Let $\chi_j$ be the characteristic function of~$V_j$.
As above, we write  $\chi_j$  instead of $\chi_j(\breve y)$.

For an arbitrary $f$ in  $\mathcal{O}(U)$, we write out the Taylor expansion in a neighborhood of the point  $\lambda_j$ and
use the bounds from Lemma~\ref{lem2}.  For some $K_j$ with $j$ for which $r_j>0$, we have
$$
\|f(\breve y)\chi_j\|\le K_j\, \|f\|_{\lambda_j,r_j,s_j},
$$
where $s_j\!:=1/(k_j-1)$  (if the order of $\lambda_j$ is equal to $1$, then
$s_j=\infty$ and $r_j\in\mathbb{N}$).  It follows readily from
the representation of holomorphic functional calculus in the form of Cauchy integral that $\sum_j\chi_j=1$. Hence $$\|f\|\le \sum_j K_j\,
\|f(\breve y)\chi_j\|,$$ and we obtain the desired assertion.
\end{proof}

We can now complete the proof of the second part of the main theorem.

\begin{proof}[Proof of Theorem~\ref{comfunenv}. Part~(B)]
 Suppose that $B$ is an Arens--Michael algebra, $\breve x,\breve y\in
B$,  the spectrum
of $\breve y$ is contained in $U$, and $[\breve x,\breve
y]=h(\breve y)$. Then, for every open  $V\subset U$  containing the spectrum, there is
a holomorphic functional calculus
$\mathcal{O}(V)\to B$  for~$\breve y$.

(1)~First, we show that
\begin{equation}\label{xbrxfbry}
[\breve
x,f(\breve y)]=h(\breve y)f'(\breve y)
\end{equation}
for every
$f\in\mathcal{O}(U)$.

Suppose first that $B$ is a Banach algebra. It can readily be seen that \eqref{xbrxfbry} holds if $f$ is a polynomial.

Let $B_0$ denote the closed subalgebra of $B$ generated by $\breve y$. Since $B_0$ is commutative, we can apply
Lemma~\ref{lem1} to the continuous derivation
$$
B_0\to B_0\!: b\mapsto [\breve x,b]\,.
$$
Thus, ${\mathop{\mathrm{Sp}}\nolimits}_{B_0} \breve y$  is finite and, therefore, so is ${\mathop{\mathrm{Sp}}\nolimits}_B \breve y$. Let $V$ be a finite union of open disks of finite radius
with pairwise disjoint closures such that ${\mathop{\mathrm{Sp}}\nolimits}_B \breve y\subset
V\subset U$.  Since the holomorphic functional calculus is
continuous, it follows that there is a compact subset $K$ of $V$ such that $|\cdot|_K$ dominates the norm on $B$
up to constant. We may assume that $K$ is a finite union of closed disks.

Since $K$ is a compact set contained in $U$, it follows that every function in $\mathcal{O}(V)$
 is continuous on $K$
and holomorphic in its interior. Moreover, the complement of $K$ is connected; therefore, by Mergelyan’s
theorem, this function can be approximated by polynomials uniformly on $K$. This implies that~\eqref{xbrxfbry} holds
for an arbitrary $f\in\mathcal{O}(V)$. The uniqueness of holomorphic functional calculus implies the equality~\eqref{xbrxfbry}  for an arbitrary $f\in\mathcal{O}(U)$.

For the general case, in which $B$ is an Arens--Michael algebra, we consider an arbitrary continuous
submultiplicative seminorm~$\|\cdot\|$ on $B$. It follows from what was proved above that the desired equality
holds in the completion with respect to~$\|\cdot\|$; in particular $$\|[\breve x,f(\breve y)]-h(\breve y)f'(\breve
y)\|=0\,.$$ Since~$\|\cdot\|$ is arbitrary, it follows that~\eqref{xbrxfbry} holds in~$B$.

(2)~Further, we show that there is a continuous homomorphism
$\nu\!:{\mathscr A}\to B$ such that $f(\breve y)=\nu\mu(f)$ for all $f\in\mathcal{O}(U)$ (here $\mu\!:\mathcal{O}(U)\to {\mathscr A}$  stands for the completion homomorphism of Proposition~\ref{compllem}). Recall that $\delta_0(f)=hf'$
for $f\in\mathcal{O}(U)$.
It follows from \eqref{xbrxfbry}  that
\begin{equation}\label{de0com}
\delta_0(f)(\breve y)=h(\breve y)f'(\breve y)=[\breve x,f(\breve
y)]\qquad (f\in\mathcal{O}(U))\,.
\end{equation}

If~$\|\cdot\|$ is a continuous submultiplicative seminorm on $B$, then
 $\|f\|_1\!:=\|f(\breve y)\|$ defines a continuous submultiplicative seminorm on  $\mathcal{O}(U)$.
It follows from \eqref{de0com} that
$$
\|\delta_0(f)\|_1=\|\delta_0(f)(\breve y)\|=\|\,[\breve x,f(\breve
y)]\,\|\le 2\,\|\breve x\|\,\| f(\breve
y)\|=2\,\|\breve x\|\,\|f\|_1\qquad (f\in\mathcal{O}(U))\,.
$$
Thus, $\|\cdot\|_1$ is $\delta_0$-stable. It follows from Proposition~\ref{cOUdernB}  that  $\|\cdot\|_1$ is dominated by the family~\eqref{prsys} of
seminorms. By Proposition~\ref{compllem}, the algebra~${\mathscr A}$  is the completion with respect to
\eqref{prsys};  hence there is a
continuous homomorphism $\nu\!:{\mathscr A}\to B$
 such that $f(\breve y)=\nu\mu(f)$ for all $f\in\mathcal{O}(U)$.

(3)~In conclusion, we show that $[\breve x,\nu(a)]= \nu\delta(a)$ for all $a\in{\mathscr A}$. Since the image of $\mu$ is dense,
it suffice to prove the equality for the case in which $a=\mu(f)$, where $f\in\mathcal{O}(U)$.
Using~\eqref{de0com}  and
Proposition~\ref{compllem}, we obtain
$$
  [\breve x,\nu(a)]=[\breve x,f(\breve y)]=\delta_0(f)(\breve y)=\nu\mu\delta_0(f)=\nu\delta\mu(f)=\nu\delta(a)\,.
$$
By Proposition~\ref{uniprOre}   there is a unique continuous homomorphism
  $\tau\!:\mathcal{O}(\mathbb{C},{\mathscr A};\delta)\to B$ such that
$\tau(x)=\breve x$ and $\nu=\tau \eta$.  It follows from the last equality that $\tau(y)=\breve y$. This completes the proof
of Theorem~\ref{comfunenv}.
\end{proof}

\subsection*{Embedding of algebras of power series in HFG  algebras}

A Fr\'echet--Arens--Michael algebra is said to be \emph{holomorphically finitely generated} or an HFG \emph{algebra} for short, if it is the quotient of the algebra of free entire functions with finitely many generators
by some closed two-sided ideal (up to topological isomorphism)
\cite[Definition~3.16, Proposition~3.20]{Pi3}.

For the sake of the completeness of our presentation, we recall that the \emph{algebra of free entire
functions} with generators $\zeta_1,\ldots,\zeta_m$ \cite{T2,T3} is the set of series
$$
\Bigl\{ a=\sum_{\alpha\in W_m} c_\alpha\zeta_\alpha\!:
\| a\|_\rho=\sum_{\alpha\in W_m} |c_\alpha|\rho^{|\alpha|}<\infty
\;\forall\rho>0\Bigr\}
$$
(with complex coefficients) equipped with the multiplication extending the concatenation operation
on the semigroup~$W_m$
of words in the alphabet $\{1,\ldots,m\}$  (for a given $\alpha\in W_m$,  the corresponding
monomial is denoted by $\zeta_\alpha$).   It can readily be seen that this algebra is a Fr\'echet--Arens--Michael
algebra.

\begin{proposition}\label{OreHFG}
For every domain $U$ in $\mathbb{C}$ and any $h\in\mathcal{O}(U)$ the algebra $\mathcal{O}(\mathbb{C},{\mathscr A};\delta)$ considered
in Theorem~\ref{comfunenv} is holomorphically finitely generated.
\end{proposition}
\begin{proof}
Although this assertion can be derived from \cite[Proposition~6.2]{Pi3},  we present a detailed proof.

Let $C$ denote the free product (or, which is the same thing, the coproduct) of  $\mathcal{O}(\mathbb{C})$ and $\mathcal{O}(U)$ in the
category of (unital) Arens--Michael algebras \cite[Sec.~4]{Pi3}. We denote the elements of $C$ corresponding
to the identity functions on $\mathbb{C}$ and $U$ by $X$ and $Y$, respectively, and consider the  closed two-sided
ideal $I$ of $C $ generated by the element $[X,Y]-h(Y)$.
Since $\mathcal{O}(\mathbb{C})$ and $\mathcal{O}(U)$ are Stein algebras and the
finiteness condition on the dimension of the embedding is satisfied, it follows that these two algebras are
HFG \cite[Theorem~3.22]{Pi3}. The property of being an HFG algebra is preserved under the formation of free
products
\cite[Corollary~4.7]{Pi3}  and passage to quotients by closed ideals
\cite[Proposition~3.18]{Pi3};  hence $C/I$
is also an HFG algebra.

We claim that $C/I$ is topologically isomorphic to
$\mathcal{O}(\mathbb{C},{\mathscr A};\delta)$. To see this, it suffices to prove that
the universal property of Theorem~\ref{comfunenv} holds. Note first that the element $Y+I$  is the image of the identity
function on $U$ under the composition $\mathcal{O}(U)\to C\to C/I$  of homomorphisms, and thus its spectrum
also belongs to $U$. Suppose further that the Arens--Michael algebra $B$ contains elements $\breve x$ and $\breve y$ such
that the spectrum
${\mathop{\mathrm{Sp}}\nolimits}_B\breve y$ is contained in $U$ and $[\breve x,\breve y]=h(\breve
y)$. Since $C$ is a free product, it follows that
the correspondence
$X\mapsto \breve x$, $ Y\mapsto \breve y$
uniquely determines a homomorphism $C\to B$,  which takes $I$
to $0$. Hence we obtain a continuous homomorphism $C/I\to B$ satisfying the desired conditions.
\end{proof}

\begin{proposition}\label{nonHFG}
For any $s\in(0,\infty]$,
 the algebra ${\mathscr A}_s$  is not HFG.
\end{proposition}
\begin{proof} (The idea of the argument below was suggested by Pirkovskii in the case of $s=\infty$; it is also
applicable for the other values of~$s$.)
Suppose that ${\mathscr A}_s$ is an HFG algebra. Since ${\mathscr A}_s$
 is commutative,
it is a Stein algebra \cite[Theorem~3.22]{Pi3}.
Since $s>0$, it is easy to see that the ideal generated by $y$ is
maximal and coincides with the Jacobson ideal. Thus,
${\mathscr A}_s$ is local. Therefore, the Gel'fand spectrum
of the algebra consists of a single point. Hence  ${\mathscr A}_s$ coincides with the algebra of germs of holomorphic
functions at this unique point of the spectrum.

On the other hand, every algebra of germs of
 holomorphic functions is a $(DF)$-space (see the case
of manifolds in the original paper by Grothendieck \cite[pp.~97--98  (Russian transl.)]{Gro} or in Mallios'
monograph \cite[pp.~136--137]{Mal}; the proof in the general case is similar). The topology on the space ${\mathscr A}_s$
is generated by a countable family of seminorms and hence is metrizable. However, every metrizable
$(DF)$-space is normable  \cite[Observation~8.3.6]{CaBo}.

If $s\in(0,\infty)$, the space ${\mathscr A}_s$  is not normable because it is isomorphic as a locally convex space to the space $\mathcal{O}(\mathbb{C})$ of entire functions, which is well known to be nonnormable (one can also apply Kolmogorov's
normability criterion \cite[II.2.1]{Scha} directly to ${\mathscr A}_s$). The space ${\mathscr A}_\infty$  cannot be normable, because it admits
no continuous norm at all. Thus, we arrive at a contradiction.
\end{proof}

Note that  ${\mathscr A}_0\cong\mathcal{O}(\mathbb{C})$ and, therefore, this is an HFG algebra.

\begin{theorem}\label{Asclsa}
 Let $s$ be a rational positive number or $\infty$. Then
${\mathscr A}_s$  is isomorphic to a closed
subalgebra of some HFG algebra.
\end{theorem}

Let us denote by $S$ the set of all positive real numbers $s$ such that ${\mathscr A}_s$
 is isomorphic to a closed
subalgebra of some HFG algebra and prove an auxiliary lemma.

\begin{lemma}\label{stspt}
If $s,t\in S$, then  $s+t\in S$.
\end{lemma}
\begin{proof}
 Suppose that  ${\mathscr A}_s$ and ${\mathscr A}_t$ are isomorphic to closed subalgebras of HFG algebras $B$ and $C$,
respectively. By the Grothendieck--Pietsch criterion
 \cite[Theorem~28.15]{MeVo}, ${\mathscr A}_s$ and ${\mathscr A}_t$
are nuclear
Fr\'echet spaces. Hence the homomorphism
${\mathscr A}_s{\mathbin{\widehat{\otimes}}} {\mathscr A}_t\to B{\mathbin{\widehat{\otimes}}} C$ is topologically injective (see, e.g., \cite[Theorem~A1.6]{EP96}). Since the class of HFG algebras is stable with respect to projective tensor products~\cite{Pi3},  it follows that ${\mathscr A}_s{\mathbin{\widehat{\otimes}}} {\mathscr A}_t$  is isomorphic to a closed subalgebra of the HFG algebra $B{\mathbin{\widehat{\otimes}}} C$. Thus, it suffices to prove that the diagonal embedding ${\mathscr A}_{s+t}\to{\mathscr A}_s{\mathbin{\widehat{\otimes}}} {\mathscr A}_t\!:y^n\mapsto  y^n\otimes  y^n$
determines a
well-defined topologically injective homomorphism of Fr\'echet algebras.

Note that ${\mathscr A}_s$   is the K\"{o}the space $\lambda(P_s)$,
corresponding to the K\"{o}the set $P_s\!:=\{r^n n!^{-s}\!:r>0\}$ (a
similar fact holds for ${\mathscr A}_t$). As noted by Pirkovskii \cite[Proposition~3.3]{Pi02},  Pietsch's results \cite{Pie63} readily
imply
$$\lambda(P_s){\mathbin{\widehat{\otimes}}}\lambda(P_t)\cong\lambda(P_s\times P_t),$$ where
$$
P_s\times
P_t=\{r^n q^m n!^{-s} m!^{-t}\!:r,q>0\}\,.
$$
Moreover, we may assume that $r=q$; thus, the natural diagonal embedding
 $\lambda(P_{s+t})\to\lambda(P_s\times P_t)$
is a
well-defined topologically injective continuous homomorphism.
\end{proof}

\begin{proof}[Proof of Theorem~\ref{Asclsa}]
Let $k\in\mathbb{N}$. By Theorem~\ref{comfunenv}, $\mathcal{O}(\mathbb{C},{\mathscr A}_{1/k};\delta)$  is the universal Arens--Michael
algebra generated by the elements $x$ and $y$ satisfying the relation  $[x,y]=y^{k+1}$.
The homomorphism
 $\eta\!:{\mathscr A}_{1/k}\to\mathcal{O}(\mathbb{C},{\mathscr A}_{1/k};\delta)$  has the form
$a\mapsto a\otimes 1$ and hence is topologically injective. Therefore, ${\mathscr A}_{1/k}$ is a closed subalgebra of
$\mathcal{O}(\mathbb{C},{\mathscr A}_{1/k};\delta)$, and the latter is an HFG algebra by Proposition~\ref{OreHFG}.Thus  $1/k\in S$.  It follows from Lemma~\ref{stspt} that all positive rational numbers belong to~$S$.

If $s=\infty$, then it suffices to apply Theorem~\ref{comfunenv}  to the relation $[x,y]=y$.
\end{proof}

For positive integer values of $s$, the assertion of Theorem~\ref{Asclsa} can be obtained also in another way, by
using the author’s results in~\cite{Ar18}.

\subsection*{Acknowledgments} The author thanks A.\,Yu.~Pirkovskii for the idea of the short proof of Proposition~\ref{nonHFG}. The author is also grateful to the referee, who noted that the case  $s=+\infty$ in Question~\ref{questionsSTsub} is not obvious, and for other
useful remarks.

\end{document}